\documentclass[12pt,a4paper]{article}
\usepackage{inputenc, amssymb, amsmath, amsthm}
\usepackage[english]{babel}
\makeatletter
\thm@style{plain}

\def\th@plain{%
  \thm@headfont{\bfseries}%
  \itshape 
  \thm@notefont{\rm}%
}
\def\thm@indent{\hspace*{\parindent}}
\makeatother

\def\({\left(}
\def\){\right)}

\newcommand{\ep}{\varepsilon}
\newcommand{\eps}{\varepsilon}

\newcommand{\be}{\begin{equation}}
\newcommand{\ee}{\end{equation}}

\let\epsilon\varepsilon
\let\phi\varphi
\let\le\leqslant
\let\ge\geqslant

\newtheorem{theorem}{Theorem}
\newtheorem{lemma}{Lemma}
\newcommand{\dokvo}{{\it Proof.} }
\newcommand\inte{\int\limits}

\begin{document}

\centerline{\bf\uppercase{A consequence of Littlewood's conditional estimates}}
\centerline{\bf\uppercase{for the Riemann zeta-function}}

\medskip

\centerline{Sergei~N.~Preobrazhenski\u i\footnote[1]{%
{\it Preobrazhenskii Sergei Nikolayevich} ---
Department of Mathematical Analysis, Faculty of Mechanics and Mathematics,
Lomonosov Moscow State University.}}

\bigskip

\hbox to \textwidth{\hfil\parbox{0.9\textwidth}{%
\small Assuming the Riemann hypothesis (RH)
and using Littlewood's conditional estimates
for the Riemann zeta-function, we provide an estimate
related to an approach of Y.~Motohashi
to the zero-free region.

\medskip

\emph{Key words}:
Riemann zeta-function, Riemann hypothesis.}\hfil}

\bigskip


{\bf 1. Introduction.} The approach of Y.~Motohashi~\cite{PreobMoto}
to the zero-free region of the Riemann zeta-function
extended by the author in~\cite{PreobPMoto} may be modified to give regions
free of large values of some products, which contain finite products $\prod_j\zeta(s_j)$.
On the Riemann hypothesis, one can obtain upper bounds for
such products for $s_j=1+t_j$ using the method of Littlewood.
To prove our result on regions free of large values we also use
an $\Omega$-theorem for $\prod_j\frac1{\zeta(s_j)}$, where
$s_j=\sigma_j+i(t_j+h_j)$
with $h_j$ lying in short intervals around $t_j$
and $\sigma_j\ge1$.
The $\Omega$-theorem depends on a version of Kronecker's theorem with an explicit upper bound.

{\bf 2. Lemmas.}
\begin{lemma} On the Riemann hypothesis,
uniformly for $\frac12<\sigma_0\le\sigma\le\frac98$ and $t\ge e^{27}$ we have
\[
\log\zeta(s)\ll\begin{cases}
\log\frac1{\sigma-1}&\text{if}\quad1+\frac1{\log\log t}\le\sigma\le\frac98,\\
\frac{(\log t)^{2-2\sigma}-1}{(1-\sigma)\log\log t}+\log\log\log t
&\text{if}\quad\sigma_0\le\sigma\le1+\frac1{\log\log t},
\end{cases}
\]
and for $\sigma>1-\frac E{\log\log t}$, $E>0$ fixed,
\be\label{PreobLTHBPbound}
\zeta(s)\ll e^{Le^{(2+\epsilon)E}}(\log\log t),
\ee
where $L=L(t)=\log\log\log\log t$
and the implied constant in the ${\ll}$ depends on $E$.
\end{lemma}

For the first estimate, see~\cite{PreobTitchHB}, Chapter XIV, \S14{.}33.
The second estimate is similar to the first and is obtained along the lines of~\cite{PreobTitchHB},
Chapter XIV, \S14{.}9. For a more precise estimate, see~\cite{PreobCarnChandee}.

\begin{lemma}\label{PreobGamma}
For $\alpha\le\sigma\le\beta$ and $t>1$ we have
\[
\Gamma(\sigma+it)=t^{\sigma+it-1/2}\exp\left(-\frac{\pi}2t-it+i\frac{\pi}2\left(\sigma-\frac12\right)
\right)\sqrt{2\pi}\left(1+O\left(\frac1t\right)\right),
\]
with the constant in the big-$O$ depending only on $\alpha$ and $\beta$.
\end{lemma}

For the proof, see e.g.~\cite{PreobKaratsVoron}, Appendix, \S3.

\begin{lemma}\label{Preobsigmaxizeta}
Let $\sigma_a(n)$, $a\in\mathbb{C}$, be the sum of $a$th powers of the divisors of $n$.
Let $\xi(d)$ be an arbitrary bounded arithmetical function
with the support in the set of square-free integers.
Then for $\sigma>1$, $T_1, T_2\in\mathbb{R}$ we have the identity
\[
\begin{split}
\sum_{n=1}^{\infty}&\sigma_{iT_1}(n)\sigma_{-iT_2}(n)
\left(\sum_{d\mid n}\xi(d)\right)n^{-s}\\
{} &=\frac{\zeta(s)\zeta(s-iT_1)\zeta(s+iT_2)\zeta(s-i(T_1-T_2))}
{\zeta(2s-i(T_1-T_2))}
\left(\xi(1)+\sum_{d=2}^\infty\xi(d)P_d(s,T_1,T_2)\right),
\end{split}
\]
where
\[
\begin{split}
&P_d(s,T_1,T_2)\\
{} &=\prod_{p\mid d}\left(1-\left(1-\frac1{p^s}\right)\left(1-\frac1{p^{s-iT_1}}\right)
\left(1-\frac1{p^{s+iT_2}}\right)\left(1-\frac1{p^{s-i(T_1-T_2)}}\right)
\left(1-\frac1{p^{2s-i(T_1-T_2)}}\right)^{-1}\right).
\end{split}
\]
\end{lemma}

\dokvo This is a version of Lemma~3 of Y.~Motohashi~\cite{PreobMoto}.
Let
\[
Z=\frac{\zeta(s)\zeta(s-iT_1)\zeta(s+iT_2)\zeta(s-i(T_1-T_2))}
{\zeta(2s-i(T_1-T_2))}.
\]
Changing the order of summation, we have
\[
\begin{split}
&\sum_{n=1}^{\infty}\sigma_{iT_1}(n)\sigma_{-iT_2}(n)
\left(\sum_{d\mid n}\xi(d)\right)n^{-s}\\
{} &=\xi(1)Z+\sum_{\substack{d\ge2,\text{$d$ square-free}\\ d=p_{d_1}\cdots p_{d_r}}}\xi(d)
\left(\sum_{k=1}^\infty
\frac{\sigma_{iT_1}(kp_{d_1}\cdots p_{d_r})\sigma_{-iT_2}(kp_{d_1}\cdots p_{d_r})}
{k^sp_{d_1}^s\cdots p_{d_r}^s}\right)\\
{} &=\xi(1)Z+\sum_{d\ge2,\text{$d$ square-free}}\xi(d)
\prod_{p\mid d}\left(\frac{\(1+p^{iT_1}\)\(1+p^{-iT_2}\)}{p^s}+
\frac{\(1-p^{i3T_1}\)\(1-p^{-i3T_2}\)}{\(1-p^{iT_1}\)\(1-p^{-iT_2}\)}\frac1{p^{2s}}+\ldots\right)\\
{} &\times\prod_{p\nmid d}\left(1+\frac{\(1+p^{iT_1}\)\(1+p^{-iT_2}\)}{p^s}+
\frac{\(1-p^{i3T_1}\)\(1-p^{-i3T_2}\)}{\(1-p^{iT_1}\)\(1-p^{-iT_2}\)}\frac1{p^{2s}}+\ldots\right)\\
{} &=\xi(1)Z+\sum_{d\ge2,\text{$d$ square-free}}\xi(d)Z
\prod_{p\mid d}\frac{
\frac{\(1+p^{iT_1}\)\(1+p^{-iT_2}\)}{p^s}+
\frac{\(1-p^{i3T_1}\)\(1-p^{-i3T_2}\)}{\(1-p^{iT_1}\)\(1-p^{-iT_2}\)}\frac1{p^{2s}}+\ldots
}
{
1+\frac{\(1+p^{iT_1}\)\(1+p^{-iT_2}\)}{p^s}+
\frac{\(1-p^{i3T_1}\)\(1-p^{-i3T_2}\)}{\(1-p^{iT_1}\)\(1-p^{-iT_2}\)}\frac1{p^{2s}}+\ldots
}.
\end{split}
\]
By an identity of Ramanujan---Wilson~\cite{PreobTitchHB}, (1{.}3{.}3),
\[
\begin{split}
\prod_{p\mid d}&\frac{
\frac{\(1+p^{iT_1}\)\(1+p^{-iT_2}\)}{p^s}+
\frac{\(1-p^{i3T_1}\)\(1-p^{-i3T_2}\)}{\(1-p^{iT_1}\)\(1-p^{-iT_2}\)}\frac1{p^{2s}}+\ldots
}
{
1+\frac{\(1+p^{iT_1}\)\(1+p^{-iT_2}\)}{p^s}+
\frac{\(1-p^{i3T_1}\)\(1-p^{-i3T_2}\)}{\(1-p^{iT_1}\)\(1-p^{-iT_2}\)}\frac1{p^{2s}}+\ldots
}\\
{} &=\prod_{p\mid d}\left(\frac{1-p^{i(T_1-T_2)-2s}}
{\left(1-p^{-s}\right)\left(1-p^{iT_1-s}\right)
\left(1-p^{-iT_2-s}\right)\left(1-p^{i(T_1-T_2)-s}\right)}-1\right)\\
{} &\times\frac
{\left(1-p^{-s}\right)\left(1-p^{iT_1-s}\right)
\left(1-p^{-iT_2-s}\right)\left(1-p^{i(T_1-T_2)-s}\right)}
{1-p^{i(T_1-T_2)-2s}}.
\end{split}
\]
This obviously ends the proof of the lemma.

\begin{lemma}\label{Preobabssigma2}
Assume the truth of the Riemann hypothesis.
Fix $E>0$. Let
\[
\exp(A\log\log T\log\log\log T)\le N\le\exp(DA\log\log T\log\log\log T),\quad T\ge e^{27},
\]
with $A=\frac{18+\eps}E$ 
and a sufficiently large positive constant $D$,
and let us put $T_1=T$, $T_2=T+H$, with $H=c(\log\log T)^{-1}$.
Then we have
\[
\begin{split}
\sum_{n\le N}&|\sigma_{iT_1}(n)|^2|\sigma_{iT_2}(n)|^2\\
{} &\mathrel{{\ll}_{A,D}}N\\
{} &\times\left((\log\log\log T)^3(\log\log T)^7|\zeta(1+iT_1)|^4|\zeta(1+iT_2)|^4\right.\\
{} &+(\log\log T)^7\zeta(1+i(T_1+H))^2\zeta(1-i(T_1-H))^2\zeta(1+i(T_2+H))^2\zeta(1-i(T_2-H))^2\\
{} &\left.+(\log\log T)^7\zeta(1+i(T_1-H))^2\zeta(1-i(T_1+H))^2\zeta(1+i(T_2-H))^2\zeta(1-i(T_2+H))^2\right)\\
{} &+O\left(N(\log\log T)^{-1}\right).
\end{split}
\]
\end{lemma}

\dokvo Let
\[
F_0(s,T_1,T_2)=\sum_{n=1}^{\infty}|\sigma_{iT_1}(n)|^2|\sigma_{iT_2}(n)|^2n^{-s}\quad
(\sigma>1).
\]
By the identity of U.~Balakrishnan~\cite{PreobBala}, we have
\[
\begin{split}
F_0(s,T_1,T_2)&=\zeta(s)^4\zeta(s+iT_1)^2\zeta(s-iT_1)^2\zeta(s+iT_2)^2\zeta(s-iT_2)^2\\
{} &\times\zeta(s+i(T_1-T_2))\zeta(s-i(T_1-T_2))\zeta(s+i(T_1+T_2))\zeta(s-i(T_1+T_2))
G(s,T_1,T_2),
\end{split}
\]
where $G(s,T_1,T_2)$ is regular and bounded for $\sigma\ge\sigma_0>1/2$,
uniformly in $T_1$, $T_2$.
The limiting case $T_1=T_2$ gives the identity of Y.~Motohashi,
which is connected with the famous nonnegative trigonometric polynomial
$3+4\cos\phi+\cos2\phi$ and the inequality of Mertens.
Littlewood's bound~\eqref{PreobLTHBPbound} and Perron's inversion formula for the height $U=N^{1+\eps}$ give
\[
\begin{split}
\sum_{n\le N}&|\sigma_{iT_1}(n)|^2|\sigma_{iT_2}(n)|^2
=\mathop{\text{\textup{Res}}}\(F_0(s,T_1,T_2)N^ss^{-1}\)_{s=1,1\pm iH}\\
{} &+O\left(\left(e^{L(T)e^{(2+\epsilon)E}}\log\log T\right)^{10}
(\log\log T)^6
N^{\eta}\log U\right)\\
{} &=\mathop{\text{\textup{Res}}}\(F_0(s,T_1,T_2)N^ss^{-1}\)_{s=1,1\pm iH}+
O\left(N(\log\log T)^{-1-\eps}\right),
\end{split}
\]
where we have put
\[
\eta=1-\frac E{\log\log T}.
\]
Also,
\[
\mathop{\text{\textup{Res}}}\(F_0(s,T_1,T_2)N^ss^{-1}\)_{s=1}\ll
N\sum_{k=0}^3|(\partial s)^k_{s=1}H(s,T_1,T_2)|(\log N)^{3-k},
\]
where
\[
\begin{split}
H(s,T_1,T_2)&=\zeta(s+iT_1)^2\zeta(s-iT_1)^2\zeta(s+iT_2)^2\zeta(s-iT_2)^2\\
{} &\times\zeta(s+i(T_1-T_2))\zeta(s-i(T_1-T_2))\zeta(s+i(T_1+T_2))\zeta(s-i(T_1+T_2)).
\end{split}
\]
By taking the logarithmic derivative, we get
\[
(\partial s)^k_{s=1}H(s,T_1,T_2)\ll H(1,T_1,T_2)(\log\log T\log\log\log T)^k.
\]
From the theorem of Littlewood and the definition of $H$ we see that
\[
\zeta(1+i(T_1-T_2))\zeta(1-i(T_1-T_2))\zeta(1+i(T_1+T_2))\zeta(1-i(T_1+T_2))\ll
(\log\log T)^4,
\]
which implies the assertion of the lemma.

\begin{lemma}\label{PreobBarbanVehov}
Let $\mu(d)$ be the M\"obius function, and let
\[
\lambda_d(z)=
\begin{cases}
\mu(d)&\text{if}\quad d<z,\\
\mu(d)\frac{\log\(z^2/d\)}{\log z}&\text{if}\quad z\le d<z^2,\\
0&\text{otherwise},
\end{cases}
\]
where $z>1$ is arbitrary. Then we have, uniformly in $N>1$ and in $z$,
\[
\sum_{n\le N}\left(\sum_{d\mid n}\lambda_d(z)\right)^2\ll\frac N{\log z}.
\]
\end{lemma}

This lemma is due to Barban---Vehov~\cite{PreobBarbanVehov}
and appears as Lemma~5 in Y.~Motohashi~\cite{PreobMoto}.
For the proof, see~\cite{PreobGraham1} and~\cite{PreobGraham2}.

\begin{lemma}\label{PreobTsanglogp}
For any large $y$, and fixed $a$, $q>1$, $(a,q)=1$,
\[
\sum_{\substack{p\le y\\ p\equiv a\pmod q}}\mathop{\text{\textup{sgn}}}
\left(\cos(2h\log p)\right)\frac{\cos(h\log p)}p=
\frac1{\phi(q)}\log\left(\min\(h^{-1},\log y\)\right)+O(1)\quad\text{\textup{for }}0<h<c.
\]
\end{lemma}

This and related estimates can be proved by using PNT in arithmetic progressions
and Stieltjes integration. A similar lemma can be found in~\cite{PreobTsang}.

{\bf 3. Proof of Theorem.}
We put
\be\label{PreobXzassum}
X=\exp(0{.}5DA\log\log T\log\log\log T),\quad z=\exp(A\log\log T\log\log\log T)
\ee
with the same $A$ and $D$ as in Lemma~\ref{Preobabssigma2},
set $\xi(d)=\lambda_d(z)$ in Lemma~\ref{Preobsigmaxizeta}
and for $T_1=T$, $T_2=T+H$ with $H=c(\log\log T)^{-1}$,
write 
\begin{gather*}
J(s,T_1,T_2)=\frac{\zeta(s)\zeta(s-iT_1)\zeta(s+iT_2)\zeta(s-i(T_1-T_2))}
{\zeta(2s-i(T_1-T_2))},\\
K(s,T_1,T_2)=\sum_{d\le z^2}\lambda_d(z)P_d(s,T_1,T_2).
\end{gather*}

\begin{theorem}
Assume the Riemann hypothesis.
Then there exists an infinite sequence of pairs of real numbers $(T_1,T_2)$,
$T_1=T$, $T_2=T+H$, with arbitrarily large values of $T$ and $H=c(\log\log T)^{-1}$,
such that
\[
|\zeta(1+iT_1)||\zeta(1+iT_2)|\ll(\log\log T)^{-2}
\]
and
\[
\begin{split}
&(\log\log T)^7|\zeta(1+iT_1)|^4|\zeta(1+iT_2)|^4\\
{}+&(\log\log T)^7\zeta(1+i(T_1+H))^2\zeta(1-i(T_1-H))^2\zeta(1+i(T_2+H))^2\zeta(1-i(T_2-H))^2\\
{}+&(\log\log T)^7\zeta(1+i(T_1-H))^2\zeta(1-i(T_1+H))^2\zeta(1+i(T_2-H))^2\zeta(1-i(T_2+H))^2\\
{}\ll&(\log\log T)^{-1}.
\end{split}
\]
Let $s_0=\sigma_0+it_0$ be a point such that
\be\label{PreobJKassumg}
|J(s_0,T_1,T_2)K(s_0,T_1,T_2)|\ge(\log\log T)^{\ep}
\ee
with arbitrarily small fixed $\ep>0$, and
\be\label{Preobsigma0t0assum}
\sigma_0=1-\frac{E_0}{\log\log T}\ge1-\frac E{\log\log T},\quad C\log\log\log T\le|t_0|\le T/2. 
\ee
Then $E_0\ge c_2(\ep)>0$.
\end{theorem}

\dokvo By Mellin's inversion formula, when $c-\sigma_0>0$,
\[
e^{-n/X}=\frac1{2\pi i}\inte_{c-i\infty}^{c+i\infty}
\Gamma(s-s_0)\frac{X^{s-s_0}}{n^{s-s_0}}\,ds.
\]
Hence for $c>1$ and $c>\sigma_0$ by Lemma~\ref{Preobsigmaxizeta} we have that
\[
\begin{split}
e^{-1/X}&+\sum_{n\ge z}\sigma_{iT_1}(n)\sigma_{-iT_2}(n)n^{-s_0}a(n)e^{-n/X}\\
{} &=\frac{X^{-s_0}}{2\pi i}\inte_{(\sigma=c)}
J(s,T_1,T_2)K(s,T_1,T_2)\Gamma(s-s_0)X^s\,ds,
\end{split}
\]
where
\[
a(n)=\sum_{d\mid n}\lambda_d(z).
\]
We now move the line of integration to the line
\[
\sigma=\eta=1-\frac E{\log\log T}.
\]
There are simple poles at $s=1$, $1+iT_1$, $1-iT_2$, $1+i(T_1-T_2)$,
but by~\eqref{Preobsigma0t0assum} and Lemma~\ref{PreobGamma} they leave residues
that are all bounded by $O\left((\log\log T)^{-2}\right)$. 
Now we consider the estimation of the integral along $\sigma=\eta$.
For the estimation of $K(s,T_1,T_2)$ we define the generating
Dirichlet series
\[
\begin{split}
M_w(s,T_1,T_2)&=1+\sum_{d=2}^{\infty}\mu(d)P_d(s,T_1,T_2)d^{-w}\\
{} &=\prod_p\left(
1-\frac1{p^w}\left(
1-\left(1-\frac1{p^s}\right)\left(1-\frac1{p^{s-iT_1}}\right)
\left(1-\frac1{p^{s+iT_2}}\right)\left(1-\frac1{p^{s-i(T_1-T_2)}}\right)\right.\right.\\
{} &\left.\left.\times\left(1-\frac1{p^{2s-i(T_1-T_2)}}\right)^{-1}
\right)
\right).
\end{split}
\]
Using a version of Perron's inversion formula, we get
\[
\frac1{1!}\sum_{d\le z^2}\mu(d)P_d(s,T_1,T_2)\log\(z^2/d\)=
\frac1{2\pi i}\inte_{c-i\infty}^{c+i\infty}M_w(s,T_1,T_2)\frac{z^{2w}}{w^2}\,dw,
\]
with $c=1-\Re s+\frac1{\log z}$,
which implies that on the line $\Re s(=\sigma)=\eta$ we have
\[
K(s,T_1,T_2)\ll z^{2(1-\eta)}(\log z)^{10}\ll\exp(2AE\log\log\log T)
(\log\log T\log\log\log T)^{10}.
\]
Thus recalling~\eqref{Preobsigma0t0assum}, \eqref{PreobJKassumg} and~\eqref{PreobXzassum} we get,
as in the proof of Lemma~\ref{Preobabssigma2}, that
\[
\begin{split}
&\left|\mathop{\text{\textup{Res}}}\(X^{-s_0}J(s,T_1,T_2)K(s,T_1,T_2)\Gamma(s-s_0)X^s\)_{s=s_0}\right.\\
{} &\left.+\frac{X^{-s_0}}{2\pi i}\inte_{(\sigma=\eta)}J(s,T_1,T_2)K(s,T_1,T_2)
\Gamma(s-s_0)X^sds-e^{-1/X}\right|\\
{} &\ge(\log\log T)^{\ep}+O\left(
\exp\left(0{.}5DA\log\log\log T(E_0-E)\right)
\frac{\log\log T}{E-E_0}\right.\\
{} &\left.\times\left(e^{L(T)e^{(2+\epsilon)E}}\log\log T\right)^{4}
(\log\log T)^{2AE+10+\eps}\right).
\end{split}
\]
Hence there is an $N$ such that $z\le N\le X^2$,
and
\[
\sum_{N\le n\le2N}|\sigma_{iT_1}(n)||\sigma_{-iT_2}(n)||a(n)|n^{-\sigma_0}\gg(\log\log T)^{-1+\ep},
\]
since the range of the summation $z\le n\le X^2$ may be divided
into the intervals $N\le n\le2N$ so that the number of the intervals
is $\ll\log X^2/z\ll\log\log T\log\log\log T$ and the sum over the entire range
must be $\gg(\log\log T)^{\ep}$.
By the Cauchy inequality and by Lemma~\ref{PreobBarbanVehov}, we get
\[
(\log\log T)^{-2+\ep}\log z\ll
\sum_{N\le n\le2N}|\sigma_{iT_1}(n)|^2|\sigma_{iT_2}(n)|^2N^{1-2\sigma_0}.
\]
Finally, by Lemma~\ref{Preobabssigma2} with $T_1=T$, $T_2=T+H$ we establish that
\[
\begin{split}
N^{2(1-\sigma_0)}&\gg
\left((\log\log\log T)^3(\log\log T)^7|\zeta(1+iT_1)|^4|\zeta(1+iT_2)|^4\right.\\
{} &+(\log\log T)^7\zeta(1+i(T_1+H))^2\zeta(1-i(T_1-H))^2\zeta(1+i(T_2+H))^2\zeta(1-i(T_2-H))^2\\
{} &\left.+(\log\log T)^7\zeta(1+i(T_1-H))^2\zeta(1-i(T_1+H))^2\zeta(1+i(T_2-H))^2\zeta(1-i(T_2+H))^2\right)\\
{} &\left.+O\left((\log\log T)^{-1}\right)\right)^{-1}(\log\log T)^{-1+\ep}.
\end{split}
\]
Next we prove existence of the infinite sequence of pairs of real numbers $(T_1,T_2)$,
claimed in the theorem.
We may choose $T_1=T$ and $T_2=T+H$ in the following way:
As in~\cite{PreobTitchHB}, Chapter VIII, \S8{.}6, for $\sigma>1$
\[
\log\frac1{|\zeta(s)|}=-\sum\frac{\cos(t\log p_n)}{p_n^{\sigma}}+O(1).
\]
Also, we have the identity
\[
\cos((t+h)\log p_n)=\cos(t\log p_n)\cos(h\log p_n)-\sin(t\log p_n)\sin(h\log p_n).
\]
So, we want to choose $t$ such that for, say, every $p_n\equiv\pm1\pmod7$ and $n\le N_2$
\[
\cos(t\log p_n)<-1+\frac1{N_2},
\]
for every $p_n\equiv\pm2\pmod7$ and $n\le N_2$
\[
\cos(t\log p_n)\begin{cases}
{}<-1+\frac1{N_2}&\text{if}\quad\cos(H\log p_n)\ge0,\\
{}>1-\frac1{N_2}&\text{if}\quad\cos(H\log p_n)<0,
\end{cases}
\]
and for every $p_n\equiv\pm3\pmod7$ and $n\le N_2$
\[
\cos(t\log p_n)\begin{cases}
{}<-1+\frac1{N_2}&\text{if}\quad\cos(2H\log p_n)\ge0,\\
{}>1-\frac1{N_2}&\text{if}\quad\cos(2H\log p_n)<0.
\end{cases}
\]
This may be done as in Lemma $\delta$ of~\cite{PreobTitchHB}, Chapter VIII, \S8{.}8.
Now existence of the sequence $(T_1,T_2)$ follows from this
and estimates as in Lemma~\ref{PreobTsanglogp}
by the Phragm\'{e}n--Lindel\"{o}f method.
Thus,
\[
\begin{split}
N^{2(1-\sigma_0)}&\gg\left((\log\log\log T)^3
(\log\log T)^{-1}
+O\left((\log\log T)^{-1}\right)\right)^{-1}\\
{} &\times(\log\log T)^{-1+\ep}.
\end{split}
\]
This ends the proof of the theorem.



\end{document}